\documentclass[12pt,a4paper]{article}
\usepackage{amsmath}
\usepackage{amssymb}
\usepackage{amsthm}
\usepackage{amsfonts}
\usepackage{latexsym}

\theoremstyle{plain}
\newtheorem{thm}{Theorem}

\newtheorem{lem}{Lemma}

\theoremstyle{definition}

\theoremstyle{remark}

\newtheorem{rem}{Remark}

\title{The asymptotic growth of equivariant sections \\of positive and big line bundles}
\author{Roberto Paoletti\footnote{\noindent{\bf Address.}
Dipartimento di Matematica "Ennio De Giorgi", Universit\'a di Lecce,
Via per Arnesano, 73100 Lecce, Italy; {\bf e-mail}:
roberto.paoletti@unile.it }}
\date{}
\begin{document}
\maketitle

\section{Introduction}

Let $X$ be a smooth complex projective $n$-fold and
$L$ a big line bundle on $X$, that is, having maximal Kodaira dimension: $\kappa (X,L)=n$. After Fujita,
the {\it volume}
$\upsilon (L)=\upsilon (X,L)$ is $$\upsilon
(X,L)=\limsup _{k\rightarrow \infty}\frac {n!}{k^n}h^0(X,{\cal
O}_X(kL)).$$ If $L$ is ample, or more generally nef and big, $\upsilon (L)=(L^n)$, the
self-intersection number of $L$. The volume of a general big line bundle has been
studied in \cite{fuj} and \cite{del}; in particular,
$\upsilon (L)$ has been given the following geometric
interpretation (Proposition 3.6 of \cite{del}): Let $(kL)^{^{[n]}}$ be the
moving self-intersection number of $kL$, that is, the number of
intersection points away from the base locus of $n$ general
divisors in the linear series $|kL|$. Then
$$\upsilon (L)=\limsup _{k\rightarrow
\infty}\frac{(kL)^{^{[n]}}}{k^n}.$$
Suppose now that a finite group $G$ acts holomorphically on $X$, and
that the action
linearizes to $L$. Let
$V_1,\ldots, V_c$ be the irreducible linear representations of $G$.
For every $k$, we have an
induced linear action of $G$ on the space of global sections
$H^0(X,L^{\otimes k})$, and therefore
an essentially unique decomposition
\begin{equation}H^0(X,L^{\otimes k})=\bigoplus _{i=1}^cH^0_i(X,L^{\otimes k}),\label{eqn:dec}\end{equation}
where each summand $H_i(X,L^{\otimes k})$ is $G$-equivariantly isomorphic to a direct sum
of copies of $V_i$. For each $i$, set $h^0_i(X,L^{\otimes k})=\dim H^0_i(X,L^{\otimes k})$ and
then define the $i$-th equivariant volume of $L$ as
\begin{equation}\label{eqn:v_i}\upsilon _i(X,L)=\limsup _{k\rightarrow \infty}\frac {n!}{k^n}h^0_i(X,L^{\otimes k}).\end{equation}
Here we shall study the volumes $\upsilon _i(L)$ and show, in particular, that

\begin{thm} \label{thm:main}
In the above situation, suppose in addition that the action of $G$ on $X$ is faithful.
Then for every $i=1,\ldots,c$ we have
\begin{eqnarray*}\label{eqn:main1}
\upsilon _i(L)=\frac {\dim (V_i)^2}{|G|}\upsilon (L).\end{eqnarray*}
\end{thm}

\noindent
For the trivial representation this has also been observed by Ein and Lazarsfeld.
Furthermore,  
if $L$ is ample Theorem \ref{thm:main} follows from algebraic results of Howe \cite{h};
more generally, if $L$ is big and numerically effective (that is, $L\cdot C\ge 0$ for every projective curve $C\subseteq X$) an algebro-geometric proof can be given applying the Riemann-Roch 
theorem on the quotient orbifold $X/G$ (I am endebted to M. Brion and J.-P. Demailly for pointing out these approaches to
me).
Here we follow however a different path, based first on an asymptotic estimate of the equivariant Szeg\"o kernels in the positive case (see below for a precise definition), and next on an equivariant version of Fujita's approximation theorem to extend the result to arbitrary big line bundles. This has the following advantages: First, as we explain below, this approach applies very naturally to the context of almost complex quantization of compact symplectic manifolds. Secondly, 
in the ample case,
it yields the lower order terms of the expansion of $h^0_i(X,L^{\otimes k})$ in decreasing powers
of $k$ (Theorem \ref{thm:main1}), in terms of the asymptotic expansion of the total Szeg\"o kernel, and under suitable assumptions on the dimension of the locus of points having non-trivial stabilizer.

As hinted above, we can  consider similarly defined invariants in the broader context of almost complex quantization.
Namely, let $(X,\omega)$ be a compact symplectic manifold with $\frac 1{2\pi}[\omega]\in H^2(X,\mathbb{Z})$,
and let $J$ be an almost complex structure on $X$ compatible with $\omega$. The pair $(\omega,J)$ fixes an Hermitian, hence a Riemannian, structure on $X$. Furthermore, by the integrality assumption on $\omega$, there exists an hermitian line bundle $L$ on $X$, having a unitary connection $\nabla _L$ of curvature $-2\pi i\omega$. Let $H(X,L^{\otimes k})\subset {\cal C}^\infty (X,L^{\otimes k})$ be the subspace of the space of smooth global sections of $L^{\otimes k}$ introduced by Guillemin and Uribe in \cite{gu}, in terms of the asymptotic spectral properties of a certain renormalized Laplacian operator; when $J$ is integral, $\omega$ a Hodge form on $X$, $L$ an ample holomorphic line bundle, and $\nabla _L$ the unique unitary connection compatible with the holomorphic structure, $H(X,L^{\otimes k})$ is the usual space of holomorphic section of $L^{\otimes k}$. The dimension of $H(X,L^{\otimes k})$ is always given, for $k\gg 0$, by the Riemann-Roch formula, and the projective embeddings defined by the linear series $|H (X,L^{\otimes k})|$ have a good asymptotic behaviour \cite{bu}, \cite{sz}.

Suppose now that the finite group $G$ acts faithfully on $X$ as a group of symplectomorphisms.
We may choose in the above a $G$-invariant compatible almost complex structure $J$, and then all the construction can be made equivariantly. Thus $G$ acts linearly on $H(X,L^{\otimes k})$, and there is a direct sum decomposition as in
(\ref{eqn:dec}): $H(X,L^{\otimes k})=\oplus _iH_i(X,L^{\otimes k})$. Setting $h(X,L^{\otimes k})=\dim H(X,L^{\otimes k})$,
$h_i(X,L^{\otimes k})=\dim H_i(X,L^{\otimes k})$, the $i$-th volume of $L$ is then defined as in (\ref{eqn:v_i}):
$\upsilon _i(L)=\limsup _{k\rightarrow +\infty}\frac {n!}{k^n}h_i(X,L^{\otimes k})$. As explained in Remark \ref{rem:acq},
the proof of the following Theorem is essentially the same as the proof of Theorem \ref{thm:ample} below.

\begin{thm}\label{thm:acq} Let $(X,\omega)$ be a compact $2n$-dimensional symplectic manifold with $\frac 1{2\pi}[\omega]$ an integral cohomology class. Choose $L$, $J$ and $\nabla$ as described above.
Suppose that the finite group $G$ acts faithfully as a group of symplectomorphisms on $X$, $J$ is $G$-invariant, and the action linearizes to $L$. Then
$$\upsilon _i(L)=\frac {\dim (V_i)^2}{|G|}\int _X\omega ^{\wedge n}.$$\end{thm}

Let us now dwell on the lower order terms of the expansion of the dimension of the covariant factors $H^0_i(X,L^{\otimes k})$, in the case where $X$ is a complex projective manifold and $L$ is ample. 
In the hypothesis of Theorem \ref{thm:main}, let $V\subset X$ be the locus of points with non-trivial stabilizer. 
Then $V$ is a union of complex submanifolds of $M$; let $c$ be its complex codimension.

\begin{thm} In the hypothesis of Theorem \ref{thm:main}, and with the above notation, assume in addition that 
$L$ is ample.
Suppose $s$ is an integer with $0\le s\le n-1$ and $c>s$. Then
$$h^0_i(X,L^{\otimes k})-\frac{\dim (V_i)^2}{|G|}h^0(X,L^{\otimes k})=o(k^{n-s})\mbox{ for every }i\,=\,1,\ldots,c.$$
If $V=\emptyset$, i.e., the action of $G$ on $X$ is free, then 
$$h^0_i(X,L^{\otimes k})=\frac{\dim (V_i)^2}{|G|}h^0(X,L^{\otimes k})\mbox{ for every }i\,=\,1,\ldots,c\mbox{ and }k\gg 0.$$
\label{thm:main1} 
\end{thm}

In particular, if $c>s$ we may compute the first $s$ terms in the asymptotic expansion of $H^0_i(X,L^{\otimes k})$ by
integrating the first $s$ terms in the asymptotic expansion of the Szeg\"o kernel of $L$ restricted to the diagonal (see below). In the
projective case the first terms of this expansion have been explicitly computed by Lu \cite{lu}.

\noindent {\bf Notation.} We shall occasionally loosely shift from multiplicative to additive notation for line bundles.
Furthermore, we shall generally identify without warning an invertible sheaf with the associated line bundle.

\noindent {\bf Acknowledgements} I am grateful to M. Brion, J.-P. Demailly and S. Zelditch for many interesting remarks.

\section{Proofs.}
To fix ideas, we focus on the complex projective case; the general almost complex case is discussed in Remark \ref{rem:acq}.
To ease the exposition, in the following by a {\it $G$-line bundle} we shall mean a line bundle on $X$ to which
the action of $G$ linearizes.

\begin{rem}
\label{rem:glb} If $D$ is a $G$-invariant divisor on $X$ (not necessarily effective), then 
${\cal O}_X(D)$ is a $G$-line bundle. 
In particular, for any line bundle $H$ on $X$, $\bigotimes _{g\in G}g^*H$ is a $G$-line bundle in a natural manner, (very) ample if so is $H$.\end{rem}

Before dealing with a general big line bundle, let us
consider the special case where $L$ is ample. 
\begin{thm} Notation being as above (with $X$ a complex projective manifold of dimension $n$), assume again that the action of $G$ on $X$ is faithful
and in addition that $L$ is an ample $G$-line bundle on $X$.
Then for every $i=1,\ldots,c$ we have
$$\upsilon _i(L)=\lim _{k\rightarrow \infty}\frac {n!}{k^n}H^0_i(X,L^{\otimes k})=
\frac {\dim (V_i)^2}{|G|}(L^n).$$ \label{thm:ample}
\end{thm}
\noindent
This is a special case of Theorem \ref{thm:main1} (the case $s=0$). We prove it separately
because it is just what is needed from the positive case to prove Theorem \ref{thm:main} in full generality, 
and furthermore its proof also
establishes, with minor modifications, Theorem \ref{thm:acq}.

\bigskip

\noindent {\it Proof.} Let $h=h_L$ be an hermitian metric on $L$ such that the curvature
form $\omega$ of the unique compatible covariant derivative $\nabla _L$ on $L$ is K\"ahler.
After averaging over $G$, we may assume that $h$, $\nabla _L$ and $\omega$ are $G$-invariant.
Thus $\omega$ is a $G$-invariant K\"ahler form, inducing a $G$-invariant volume form $dx$ on $X$.

Let $L^*=L^{-1}$ be the dual line bundle, with the induced heritian structure and connection, and consider the unit disc bundle $L^*\supset\mathbb{S}\stackrel{\pi}{\rightarrow}X$. Let
$i\alpha\in \Omega ^1_\mathbb{S}(i\mathbb{R})$ be the connection form.
Then $dp=:\frac 1{2\pi}\pi ^*(dx)\wedge \alpha$ is a $G$-invariant volume form on $\mathbb{S}$.
For every integer $k\ge 0$, denote by
$$P_k:L^2(\mathbb{S})\rightarrow \tilde H_k(\mathbb{S})\cong H^0(X,L^{\otimes k})$$
the orthogonal projection onto the space of boundary values of holomorphic functions in the $k$-th isotype with respect to the $S^1$-action. For each $i=1,\ldots,c$, let
$$P_{k,i}:L^2(\mathbb{S})\rightarrow \tilde H_{k,i}(\mathbb{S})\cong H^0_i(X,L^{\otimes k})$$
denote the orthogonal component onto the $i$-th isotype of $\tilde H_k(\mathbb{S})$ with respect to the $G$-action.
Also, let $\tilde \Pi _k,\,\tilde \Pi_{k,i}\in {\cal C}^\infty (\mathbb{S}\times \mathbb{S})$ be the Schwartz kernels of $\Pi _k$ and $\Pi _{k,i}$, respectively. Clearly, $\tilde \Pi _k=\sum _i\tilde \Pi _{k,i}$ and
$$H^0_i(X,A^{\otimes k})=\int _{\mathbb{S}}\tilde \Pi _{k,i}(p,p)\,dp.$$
Now $\tilde \Pi _k(p,p)$ and $\tilde \Pi _{k,i}(p,p)$ descend to positive functions $\nu _k(x)$ and $\nu _{k,i}(x)$ on $X$ \cite{sz}, \cite{p}.
Hence,
\begin{eqnarray*}
\upsilon _i(L)=\limsup _{k\rightarrow \infty}\frac {n!}{k^n}\int _X\nu _{k,i}(x)\,dx.
\end{eqnarray*}
As in \cite{p}, we decompose $ \nu _{k,i}(x)$ as the sum of two terms, the first being a multiple of
$ \nu _k(x)$ and the second growing at most as $k^{(n-1)/2}$. More precisely, if $G_x\subseteq G$ is the stabilizer subgroup of $x\in X$ and $p\in \mathbb{S}$ is any point lying over $x$, we have
\begin{eqnarray}
\label{eqn:G_x} \nu _{k,i}(x)=\frac {\dim (V_i)}{|G|}\cdot (\alpha
_x^k,\chi _i)_{{}_{G_x}}\cdot \nu _k(x)+\frac {\dim (V_i)}{|G|}
\sum _{g\not\in G_x}\overline \chi _i(g)\tilde \Pi
_N(g^{-1}p,p).\end{eqnarray} Here notation is as follows: $\chi
_i:G\rightarrow \mathbb{C}$ is the character of the irreducible
representation $V_i$, $\alpha _x:G_x\rightarrow S^1\subset
\mathbb{C}^*$ is the unitary character describing the action of
$G_x$ on $L(x)$ (the fibre of $L$ over $x$), and
$(h,k)_{{}_{G_x}}=\sum _{g\in G_x}f(g)\cdot \overline k(g)$ is
the $L^2$-hermitian product with respect to the counting measure
on $G_x$. Furthermore, there exists $a>0$ such that, setting $d_x=\min \{{\rm dist}(x,gx):g\not
\in G_x\}$, the latter term is bounded above by $Ck^ne^{-a\sqrt kd_x}$ (see section 6 of \cite{christ}).

If, in particular, $G_x=\{e\}$, where $e\in G$ is the unit, the former term is
$(\dim (V_i)^2/|G|)\cdot \nu _k(x)$.
Let us now recall the following simple useful fact.

\begin{lem}
\label{lem:useful}
Suppose that the finite group $G$ acts faithfully and holomorphically on
the connected projective manifold $X$.
Then there is a nonempty
Zariski dense open subset $U\subseteq X$ such that $G_x=\{e\}$ for every $x\in U$.\end{lem}

Set $Z=X\setminus U$; thus $Z$ is a proper algebraic subvariety of $X$ of codimension, say, $c$.
For $\epsilon >0$ let $V_\epsilon \subseteq X$ be the $\epsilon$-neighbourhood of $Z$
in the geodesic distance associated to $\omega$. Then $V_\epsilon$ has volume
$\le C\epsilon ^{2c}$, where $C$ is a constant. On the other hand, by the above and the asymptotic expansion of $\nu _k(x)$ in Theorem 1 of \cite{z},
 $n!\,k^{-n}\nu _{k,i}(x)$
is in any event a bounded function. Therefore,
\begin{eqnarray*}\label{estimate:1}
\left |\frac {n!}{k^n}\int _X\nu _{k,i}(x)\,dx-
\frac {n!}{k^n}
\int _{X\setminus V_\epsilon}\nu _{k,i}(x)\,dx\right |=\left |\frac {n!}{k^n}\int _{V_\epsilon}\nu _{k,i}(x)\,dx\right |
\le C\epsilon ^{2c}.\end{eqnarray*}
There exists $\delta =\delta _\epsilon>0$ such that ${\rm dist}(x,gx)>\delta$ if
$x\in X\setminus V_\epsilon$ and $g\neq e$. Hence
\begin{eqnarray*}
\frac {n!}{k^n}\left |\int _{X\setminus V_\epsilon}\nu _{k,i}(x)\,dx-\frac {\dim (V_i)^2}{|G|}
\int _{X\setminus V_\epsilon}\nu _k(x)\,dx\right |\le Ce^{-\sqrt k\delta }.
\end{eqnarray*}
Summing up,
\begin{eqnarray*}
\left |\frac {n!}{k^n}\int _X\nu _{k,i}(x)\,dx-\frac {n!}{k^n}\cdot \frac {\dim (V_i)^2}{|G|}\int _X\nu _k(x)\,dx\right |\le C(\epsilon ^{2c}+
e^{-\sqrt k\delta }).\end{eqnarray*}
Since
$ n!\,k^{-n}\int _X\nu _k(x)\,dx\rightarrow {\rm vol}(L)=(L^n)$, the statement follows by taking
$k\gg  1/\epsilon \gg 0.$

\begin{rem} \label{rem:acq}
The asymptotic expansions used in the proof also hold in almost complex quantization \cite{bu}, \cite{sz}, \cite{p}, but the off-diagonal estimate on the $k$-th Fourier coefficient of the Szeg\"o kernel from \cite{christ} has been proved only in the complex projective case. However, in the more general almost complex case we still have the estimate
$$\left |\Pi _k(x,y)\right |\le C\nu _k(x)e^{-k{\rm dist}(x,y)^2/2}+O(k^{(n-1)/2}),$$
from \cite{bu} and \cite{sz}, which is still enough to prove the theorem.
Furthermore, since $G$ preserves the Riemannian structure on $X$ associated to $\omega$ and $J$,
in place of Lemma \ref{lem:useful} we may as well use Theorem 8.1 on page 213 of \cite{st}:
the set of all $x\in X$ with non-trivial stabilizer is a finite collection of proper submanifolds (the action being faithful). The same argument, with minor changes, thus also proves Theorem \ref{thm:acq}.
\end{rem}

\begin{lem} \label{lem:tech}
In the same hypothesis, let $A$ and $B$ be $G$-line bundles on
$X$, with $A$ ample. Then $H^0_i(X,A^{\otimes m}\otimes B)\neq 0$
for every $i=1,\ldots,c$ and $m\gg 0$.\end{lem}

\noindent {\it Proof.} Let $H$ be any ample $G$-line bundle on
$X$. Then, perhaps after replacing $H$ by $H^{\otimes k|G|}$ for
some fixed $k\gg 0$, we may assume that the linear series $\left
|H^0(X,H)^G\right |$ (corresponding to the subspace of
$G$-invariant sections of $H$) is base point free. In fact, by
\cite{p}, for any $i=1,\ldots,c$ and $k\gg 0$ the base locus of
$\left |H^0_i(X,H^{\otimes k})\right |$ is contained in the locus
$\{x\in X:(\alpha _x^k,\chi _i)_{{}_{G_x}}=0\}$ (notation here is
as in the proof of Theorem \ref{thm:ample}, with $L=H$). If $|G|$
divides $k$ then $(\alpha _x^k,\chi _i)=\sum _{g\in G_x}\overline
\chi _i(g)$; if $V_i$ is the trivial representation this is
$|G_x|\neq 0$.

Let then $V_1,\ldots, V_{n-1}\in \left |H^0(X,H)^G\right |$ be
general divisors; their complete intersection is a smooth
$G$-invariant curve $C\subseteq X$. As $A$ is ample, for $m\gg 0$
the restriction $$H^0(X,A^{\otimes m}\otimes B)\longrightarrow
H^0(C,A^{\otimes m}\otimes B\otimes {\cal O}_C)$$ is a surjective
$G$-equivariant linear map. Hence for every $i=1,\ldots,c$ we have
surjective maps $$H^0_i(X,A^{\otimes m}\otimes B)\longrightarrow
H^0_i(C,A^{\otimes m}\otimes B\otimes {\cal O}_C).$$ On the other
hand, $A^{\otimes m}\otimes B\otimes {\cal O}_C=W_m^{\otimes
d_m}$, where $W_m$ is a line bundle of degree one on $C$ (hence
ample) and $d_m=m(A\cdot C)+(B\cdot C)$. In other words, for
$m\gg 0$ and every $i$ there are surjections $$H^0_i(X,A^{\otimes
m}\otimes B)\longrightarrow H^0_i(C,W_m^{\otimes d_m}).$$ Since
${\rm Pic}^1(C)$ is compact, for any $\epsilon >0$ there is a
uniform estimate $$H^0_i(C,W_m^{\otimes d_m})\ge d_m\left (\frac
{\dim (V_i)^2}{|G|}-\epsilon \right )$$ for $m\ge m_\epsilon$.
This completes the proof.

\begin{lem} \label{lem:fuj}
Let $X$ be a smooth complex projective $n$-fold, $G$ a finite group acting
holomorphically on $X$ and $L$ a big $G$-line bundle on $X$.
Then for any $G$-line bundle $H$ on $X$ there exists $m_0\in \mathbb{N}$ such that for every
integer $m\ge m_0$ and every $i=1,\ldots,c$ we have
$$\upsilon _i\big (L^{\otimes m}\otimes H^{-1}\big )\ge m^n\big (\upsilon _i(L)-\epsilon\big ).$$
\end{lem}

\noindent {\it Proof.} This extends to our setting Lemma 3.5 in
\cite{del}, and the proof only requires some slight modifications
to the argument given there. Fix $i\in \{1,\ldots,c\}$ and
$\epsilon
>0$.  By definition, there is a sequence $k_\nu \uparrow +\infty$
such that $$h^0_i(X,L^{\otimes k_\nu})\ge \frac {k_\nu
^n}{n!}\left (\upsilon _i(L)-\frac \epsilon 2 \right ).$$ Fix
$m\gg 0$ and set $\ell _\nu=\left [\frac {k_\nu}m\right ]$,
$r_\nu=k_\nu-\ell _\nu m$ so that, in additive notation, $$\ell _\nu (mL-H)=k_\nu
L-(r_\nu L+\ell _\nu H).$$

After replacing $H$ by $H\otimes E$ for a suitably positive
$G$-line bundle $E$ on $X$, we may suppose that $H$ is a very
ample $G$-line bundle. Furthermore, as in the proof of Lemma
\ref{lem:tech}, perhaps after replacing $H$ by a suitably large
tensor power of $H^{\otimes |G|}$, we may also assume that $\left
|H^0(X,H)^G\right |$ is base point free. Choose a smooth divisor
$D\in \left |H^0(X,H)^G\right |$. Since $D$ is a $G$-invariant
submanifold of $X$, perhaps after a change of linearization
we may also suppose that the bundle action of $G$ on $H$ is the
natural action on ${\cal O}_X(D)$ induced by the action on $k(X)$. Thus  
for every $G$-line bundle $A$ on $X$ we have a
$G$-equivariant isomorphism $$H^0(X,A(-D))=H^0(X,A\otimes {\cal I}_D)\cong H^0(X,A\otimes H^{-1}),$$
where ${\cal I}_Z\subseteq {\cal O}_X$ denotes the ideal sheaf of a closed subscheme $Z\subseteq X$. 
Furthermore, for every $i=1,\ldots,c$ the short exact sequence of sheaves $$0\longrightarrow
A(-D)\longrightarrow A\longrightarrow A\otimes {\cal O}_D\longrightarrow 0$$ induces an exact sequence
\begin{equation}\label{eqn:H^0_i}0\rightarrow H^0_i(X,A(-D))\rightarrow
H^0_i(X,A)\rightarrow H^0_i(D,A\otimes {\cal O}_D).\end{equation} Finally, by the statement of Lemma
\ref{lem:tech} we may also assume that $$\left |H^0(X,H^{\otimes
b}\otimes L^{-r})^G\right |\neq \emptyset,$$ for all integers
$0\le r\le m-1$ and $b\ge 1$. If $\sigma \in H^0(X,H^{\otimes
m}\otimes L^{-r_\nu})^G$ is non-zero, in additive notation tensor product by $\sigma$
determines injections $$H^0_i\big (X,L^{\otimes k_\nu}(-(\ell
_\nu +m)D)\big )\hookrightarrow H^0_i\big (X,L^{\otimes m\ell _\nu}(-\ell _\nu D)\big )$$ for every $\nu$; therefore $$H_i\big (X,{\cal
O}_X(\ell _\nu(mL-H))\big )\ge H^0_i\big (X,{\cal O}_X(k_\nu L-
(\ell _\nu +m)H)\big ).$$ Now we set $A=L^{\otimes k_\nu }(-jD)$ in
(\ref{eqn:H^0_i}), and proceed inductively as in {\it loc. cit.},
Lemma 3.5. More precisely, for any $s>0$ the exact sequences
\begin{eqnarray*}
0\rightarrow L^{\otimes k_\nu}(-(j+1)D)\rightarrow L^{\otimes k_\nu}(-jD)\rightarrow  L^{\otimes k_\nu}(-jD)\otimes {\cal O}_D\rightarrow
0\end{eqnarray*} for $0\le j<s$ imply
\begin{eqnarray*} H^0_i\big (X,L^{\otimes k_\nu}(-sD)\big )\ge
H_i\big (X, L^{\otimes k_\nu}\big
)-\sum _{0\le j<s}H^0_i\big (D,L^{\otimes k_\nu}(-jD)\otimes {\cal O}_D\big )\\
\ge H^0_i\big (X,L^{\otimes k_\nu}\big )-s\,H^0(D,L^{\otimes k_\nu}\otimes {\cal O}_D) \ge \frac {k_\nu ^n}{n!}\left
(\upsilon _i(L)-\frac \epsilon 2\right
)-sCk^{n-1}_\nu.\end{eqnarray*}The statement follows as in
\cite{del} by letting $\ell _\nu \gg m\gg 1$.

\bigskip
In particular, if $L$ is any big $G$-line bundle on $X$, for any
$\epsilon
>0$ there is $m_0\in \mathbb{N}$ such that \begin{equation}\upsilon _i(L)\ge
m^{-n}{\upsilon _i(mL)}\ge \upsilon
_i(L)-\epsilon\label{eqn:asympt}\end{equation} for every integer
$m\ge m_0$. Unless $V_i$ is the trivial representation, the
spaces $H^0_i(X,L^{\otimes k})$ do not form a graded linear
series. Therefore homogeneity of $\upsilon _i$ does not follow
directly from Lemma 3.4 of \cite{eln}.

\begin{lem} Let $L$ be a nef and big $G$-line bundle on $X$.
Then $$\upsilon _i(L)=\frac {\dim (V_i)^2}{|G|}(L^n).$$
\label{lem:nefand big}
\end{lem}

\noindent {\it Proof.} By the definition of $\upsilon _i$, $\sum
_i\upsilon _i(L)\ge \upsilon (L)$. Fix rational numbers $\epsilon
,\delta >0$ and let $A$ be an ample $G$-line bundle on $X$. Let
$r\gg 0$ be an integer such that $r\delta \in \mathbb{N}$. By
choosing $r$ sufficiently large and divisible, we may assume that
there exists a $G$-invariant non-zero section $\sigma \in H^0\big
(X,{\cal O}_X(r\delta A)\big )^G$. Tensor product by $\sigma$
determines  for every $i$ injective maps $$H^0_i(X,{\cal
O}_X(rL))\hookrightarrow H^0_i(X,{\cal O}_X(r(L+\delta A)).$$
Since $L+\delta A$ is ample, for $r\gg 0$ we have
\begin{eqnarray*}\frac{\dim (V_i)^2}{|G|}r^n\big ((L+\delta A)^n\big )=
\upsilon _i\big (r(L+\delta A)\big )\ge \upsilon _i(rL)\ge
r^n(\upsilon _i(L)-\epsilon),\end{eqnarray*} and taking $\epsilon$
and $\delta$ arbitrarily small we conclude that $$ \frac{\dim
(V_i)^2}{|G|}(L^n)\ge \upsilon _i(L)$$
for every $i$. Thus, since
$\sum _i\dim (V_i)^2=|G|$,
\begin{eqnarray*}\upsilon (L)=(L^n)=\sum _i\frac{\dim
(V_i)^2}{|G|}(L^n)\ge \sum _i\upsilon _i(L)\ge \upsilon
(L).\end{eqnarray*} This implies the statement.

\begin{lem} Let $L$ be any big $G$-line bundle on $X$.
Then $\upsilon _i(L)>0$ for every $i=1,\ldots,c$.
\label{lem:>0}\end{lem}

\noindent {\it Proof.} Fix $m\gg 0$ with $\dim \phi _m(X)=n$,
where $$\phi _m:X-\rightarrow \mathbb{P}H^0(X,L^{\otimes m})^*$$
is the rational map associated to the linear series $|L^{\otimes
m}|$. Let $\psi :X'\rightarrow X$ be a $G$-equivariant resolution
of singularities of $|L^{\otimes m}|$ \cite{aw}. Then $$|\psi
^*(L^{\otimes m})|=|M|+F,$$ where $F$ is the fixed divisor of
$|\psi ^*(L^{\otimes m})|$ and $M$ is a base point free (hence
nef) big $G$-line bundle on $X'$. For every $k$ we have
$G$-equivariant injective maps $$H^0\big (X',M^{\otimes k}\big
)\longrightarrow H^0\big (X',\psi ^*(L^{\otimes mk})\big ) \cong
H^0(X,L^{\otimes mk}).$$ Therefore, $\upsilon _i(L)\ge
m^{-n}\upsilon _i(L^{\otimes m})\ge m^{-n}\upsilon _i(M)$ for
every $i$. The statement then follows from the nef and big case
of Lemma \ref{lem:nefand big}.

\bigskip

\begin{lem} Let $L$ be a big $G$-line bundle on $X$. Then for every $i=1,\ldots,c$
and every $m\gg 0$ there exists $D\in |H^0_i(X,L^{\otimes m})|$
which can be written as $D=A+E$, where $A\subset X$ is a
$G$-invariant ample divisor and $E\subset X$ is an effective
divisor such that ${\cal O}_X(E)$ is an effective $G$-line bundle
with $E\in |H^0_i(X,{\cal O}_X(E))|$.\label{lem:zariski}\end{lem}

\bigskip

\noindent {\it Proof.} Fix a $G$-invariant very ample smooth
divisor $A\subseteq X$ and let $\sigma \in H^0(X,{\cal O}_X(A))^G$
be a section with $A={\rm div}(\sigma)$. Consider the short exact
sequence $$0\longrightarrow H^0_i\big (X,L^{\otimes m}(-A)\big
)\stackrel{\otimes \sigma}{\longrightarrow }H^0_i(X,L^{\otimes m}
)\longrightarrow H^0_i(A,L^{\otimes m}\otimes {\cal O}_A).$$ By
Lemma \ref{lem:>0}, $H^0_i(X,L^{\otimes m})=O(m^n)$; since
$H^0_i(A,L^{\otimes m}\otimes {\cal O}_A)\le Cm^{n-1}$, we
conclude that $H^0_i(X,L^{\otimes m}(-A))\neq 0$ for $m\gg 0$.

\bigskip
By taking $V_i$ to be the trivial representation, we see in
particular that there exists $D\in |L|$ of the form $D=A+E$,
where $A$ and $E$ are $G$-invariant $\mathbb{Q}$-divisors, with
$A$ ample and $E$ effective. If $m\in \mathbb{N}$ is such that
$mA$ and $mE$ are integral, we obtain $G$-invariant injections
$H^0(X,{\cal O}_X(mkA))\rightarrow H^0(X,L^{\otimes km})$ for
every $k$, whence \begin{eqnarray*}\upsilon _i(L)\ge
m^{-n}\upsilon _i(L^{\otimes m})\ge m^{-n}\upsilon _i\big
({\cal O}_X(mA)\big )\\=m^{-n}\frac {\dim (V_i)^2}{|G|}\big
((mA)^n\big )=\frac {\dim (V_i)^2}{|G|}(A^n).\end{eqnarray*}
\noindent
Similarly, of course, $\upsilon (L)\ge (A^n)$.

Theorem \ref{thm:main} is now a consequence of the following
equivariant version of a Theorem of Fujita \cite{del}, \cite{fuj}.

\begin{thm} Let $X$ be a smooth projective $n$-fold, $G$ a finite group acting
holomorphically and faithfully on $X$. Let $L$ be a big $G$-line
bundle on $X$. Fix $\epsilon >0$. Then there exists a
$G$-equivariant birational modification (depending on $\epsilon$)
$$\mu :X'\longrightarrow X$$ and a decomposition $\mu^*(L)\equiv
E+A$, where $E$ and $A$ are $G$-invariant $\mathbb{Q}$-divisor on
$X'$, with $E$ effective and $A$ ample, such that $(A^n)\ge
\upsilon (X,L)-\epsilon$ and $\frac {\dim (V_i)^2}{|G|}(A^n)\ge
\upsilon _i(X,L)-\epsilon$.
\end{thm}

\noindent {\it Proof.} By Fujita's Theorem, for every $\epsilon >0$
there exist a birational modification
$\mu :X'\rightarrow X$ and a decomposition $\mu ^*(L)\equiv A+E$, 
with $\mathbb{Q}$-divisors $A$ and $E$, ample and effective respectively, 
such that $(A^n)\ge \upsilon (X,L)-\epsilon$. 
Thus, in order to ensure the second inequality for every $i$
we need only give an equivariant version of the proof in \cite{del}, Theorem 3.2.

To this end, it suffices to produce for every
$\epsilon >0$ a birational modification $\mu :X'\rightarrow X$ and a decomposition
$\mu ^*(L)\equiv A+E$, where $A$ and
$E$ are $G$-invariant $\mathbb{Q}$-divisors on $X'$, $A$ is big and nef
and $E$ is effective, satisfying the stated numerical conditions.
In fact, by Lemma \ref{lem:zariski}, $A\equiv A'+D$, where $A'$
and $D$ are ample and effective $G$-invariant divisors,
respectively. Therefore, for any $\delta \in \mathbb{Q}_+$ we
have $A+E\equiv A''+F$, where $A''=(1-\delta)A+\delta A'$ is
ample, $F=E+\delta D$ is effective and $((A'')^n)$ approximates
$(A^n)$ as closely as desired.

Let $B$ be a $G$-line bundle on $X$, so positive that
$R=:K_X\otimes B^{\otimes (n+1)}$ is very ample and
$H(X,R)^G\neq \{0\}$ (Lemma \ref{lem:tech}). Fix a non-zero
section $\sigma \in H(X,R)^G$ and set $M_m=L^{\otimes m}\otimes
R^{-1}$. Then for $m$ sufficiently large $M_m$ is a big $G$-line
bundle and $\upsilon _i(M_m)\ge m^n(\upsilon _i(L)-\epsilon)$
(Lemma \ref{lem:fuj}). Tensoring with $\sigma ^{\otimes \ell}$
determines for every $\ell \ge 1$ and $i=1,\ldots,c$ injective
linear maps $$H^0_i\left (X,{\cal O}_X(M_m^{\otimes \ell })\right
)\longrightarrow H^0_i\big (X,{\cal O}_X(L^{\otimes \ell m})\big ),$$ whence
$\upsilon _i(L^{\otimes m})\ge \upsilon (M_m)$. Summing up,
$$m^n\upsilon _i(L)\ge \upsilon _i(L^{\otimes m})\ge \upsilon
_i(M_m)\ge m^n\big (\upsilon _i(L)-\epsilon \big ).$$ Let us now
consider the asymptotic multiplier ideal \cite{del}
$${\cal J}={\cal J}(X,||M_m||).$$ Then ${\cal J}={\cal J}(\frac
1k|kM_m|)={\cal J}(\frac 1k\mathfrak b_{kM_m})$ for $k\gg 0$, where
$\mathfrak b_{kM_m}\subset {\cal O}_X$ is the base ideal of the
linear series $\left |M_m^{\otimes k}\right |$. As $M_m$ is a
$G$-line bundle, $\cal J$ is a $G$-invariant ideal sheaf. Let $\mu
:X'\rightarrow X$ be a $G$-equivariant log-resolution of ${\cal
J}$ \cite{aw}, so that $\mu ^*{\cal J}={\cal O}_{X'}(-E_m)$ for
some $G$-invariant effective divisor $E_m$ on $X'$.

Since $$L^{\otimes m}\otimes {\cal J}(||M_m||)=M_m\otimes
K_X\otimes B^{\otimes (n+1)}\otimes {\cal J}(||M_m||),$$ is
globally generated by Theorem 1.8 of \cite{del}, so is the
$G$-line bundle $A_m=:\mu ^*(L^{\otimes m})(-E_m)$ on $X'$.

Since all the sheaves involved are $G$-sheaves and $\sigma$ and $E_m$ are
$G$-invariant, using Theorem 1.8 (iii) of {\it loc. cit.},
subadditivity and tensor product by $\sigma ^{\otimes \ell}$ we
have a chain of $G$-equivariant inclusions
\begin{eqnarray*}H^0\big (X,M_m^{\otimes \ell}\big )\cong H^0\big (X,M_m^{\otimes \ell}\otimes
{\cal J}(||M^{\otimes \ell}||)\big ) \subseteq H^0(X,M_m^{\otimes
\ell }\otimes {\cal J}^\ell)\\ \subseteq H^0(X,L^{\otimes
m\ell}\otimes {\cal J}^{\ell})\subseteq H^0(X',L^{\otimes
m\ell}(-\ell E_m))=H^0(X',A_m^{\otimes \ell}).\end{eqnarray*}

Thus $A_m$ is a nef and big $G$-line bundle on $X'$; we have
$$(A_m^n)=\upsilon (X',A_m)\ge \upsilon (X,M_m)\ge m^n\big (\upsilon (L)-\epsilon\big )$$
and 
$$\frac {\dim (V_i)^2}{|G|}(A_m^n)=\upsilon _i(X',A_m)\ge \upsilon _i(X,M_m)\ge m^n\big (\upsilon _i(L)-\epsilon\big )$$
for every $i=1,\ldots,c$.
Now we are done: we need only choose some $D_m\in \left |H(X',A_m)^G\right |$ and set 
$A=\frac 1m D_m$, $E=\frac 1mE_m$.

\section{Proof of Theorem 3.}

\noindent
By theorem 8.1 on page 213 of \cite{st}, $V$ is a union of submanifolds. 
If $H\subset X$ is a subgroup and $V_H\subset X$ is the submanifold of the points fixed by $H$, around any $p\in V_H$
there are local coordinates in terms of which every $g\in H$ is a linear transformation, and therefore $V_H$ is
a linear subspace. Therefore, $H$ acts freely on the unit sphere bundle of the normal bundle of $V_H$, and this
implies that there is $a>0$ such that ${\rm dist}(gx,x)\ge a\, {\rm dist}(x,V_H)$, for any $x$ sufficiently close to
$V_H$ and every $g\in H\setminus \{e\}$.
It follows, in the notation of the Theorem, that there exists $a>0$ such that if $x\not\in V_\epsilon$ for sufficiently
small $\epsilon >0$, then ${\rm dist}(gx,x)\ge a\epsilon$ for $g\neq e$ ($e$ is the neutral element of $G$).
The constants involved in the coming estimates will be allowed to vary from line to line without mention.

In view of the above and the off-diagonal estimate on the Szeg\"o kernel discussed in section 6 of \cite{christ},
if $x\not\in V_\epsilon$ and $g\neq e$ then $$\left |\Pi _k(gx,x)\right |\le  Ck^ne^{-a\sqrt k\epsilon}.$$
Arguing as in the proof of Theorem \ref{thm:main}, we have the following estimates:
\begin{eqnarray*}
\left |h^0_i(X,L^{\otimes k})-\frac{\dim (V_i)^2}{|G|}h^0(X,L^{\otimes k})\right |=
\left |\int _X\left (\nu _{k,i}(x)-
\frac{\dim (V_i)^2}{|G|}\nu _k(x)\right )dx\right |\le \\
\int _{V_\epsilon}\left | \nu _{k,i}(x)-
\frac{\dim (V_i)^2}{|G|}\nu _k(x)\right |dx+
\int _{X\setminus V_\epsilon}\left |\nu _{k,i}(x)-
\frac{\dim (V_i)^2}{|G|}\nu _k(x)\right |dx\le \\
Ck^n(\epsilon ^{2c}+e^{-a\sqrt k\epsilon}).\end{eqnarray*}
\noindent
If now $\alpha \in (\frac sc,1)$ and $\epsilon =k^{-\alpha /2}$, we have
$k^n(\epsilon ^{2c}+e^{-a\sqrt k\epsilon})=o(k^{n-s})$ 
as $k\rightarrow \infty$.

\end{document}